\input amstex 
\input amsppt.sty 
\magnification=\magstep1
\parindent=1em
\baselineskip 15pt
\hsize=12.3cm
\vsize=18.5cm
\def\nmb#1#2{#2}         
\def\cit#1#2{\ifx#1!\cite{#2}\else#2\fi} 
\def\totoc{}             
\def\idx{}               
\def\ign#1{}             

\def\colon{{:}\;}
\redefine\o{\circ}

\define\al{\alpha}

\define\et{\eta}

\define\si{\sigma}

\define\ph{\varphi}

\define\om{\omega}

\define\La{\Lambda}

\define\Ph{\Phi}

\define\Om{\Omega}
\predefine\ii{\i}
\redefine\i{^{-1}}
\define\row#1#2#3{#1_{#2},\ldots,#1_{#3}}
\define\x{\times}
\define\M{\Cal M}

\define\Vol{\operatorname{vol}}

\define\tr{\operatorname{tr}}
\define\dd#1{\frac{\partial}{\partial#1}}

\def\today{\ifcase\month\or
 January\or February\or March\or April\or May\or June\or
 July\or August\or September\or October\or November\or December\fi
 \space\number\day, \number\year}
\topmatter
\title Geodesics on spaces \\
of almost hermitian structures \endtitle
\author  Olga Gil-Medrano\footnote{Research 
partially supported by the CICYT grant n. PB91-0324,
\hfill\hfill\hfill} \\
Peter W. Michor  \footnote{Supported by Project P 7724 PHY
of `Fonds zur F\"orderung der wissenschaftlichen Forschung'.
\hfill\hfill\hfill}
\endauthor
\leftheadtext{O\. Gil-Medrano, P\. W\. Michor}
\affil
Departamento de Geometr\'\ii a y Topolog\'\ii a\\
Universidad de Valencia, Spain.\\
\\
Institut f\"ur Mathematik, Universit\"at Wien, Austria\\
Erwin Schr\"odinger Institute for Mathematical Physics.
\endaffil
\address{O\. Gil-Medrano:
Departamento de Geometr\'\ii a y Topolog\'\ii a,
Facultad de Ma\-te\-m\'a\-ti\-cas,
Universidad de Valencia,
46100 Burjassot,
Valencia, Spain}\endaddress
 
\address{P\. W\. Michor:
Institut f\"ur Mathematik, Universit\"at Wien,
Strudlhofgasse 4, A-1090 Wien, Austria; 
`Erwin Schr\"odinger International Institute of Mathematical 
Physics', Pasteurgasse 6/7, A-1090 Wien, Austria.
}\endaddress
\email{michor\@pap.univie.ac.at} \endemail
\date{\today}\enddate
\abstract{A natural metric on the space of all almost hermitian 
structures on a given manifold is investigated.}\endabstract
\endtopmatter
\document
 
\head Table of contents \endhead
\noindent 0. Introduction \leaders \hbox to 1em{\hss .\hss }\hfill {\eightrm 1}\par 
\noindent 1. Almost hermitian structures \leaders \hbox to 1em{\hss .\hss }\hfill {\eightrm 2}\par 
\noindent 2. The geodesic equation in $\Cal H$ \leaders \hbox to 1em{\hss .\hss }\hfill {\eightrm 5}\par 
\noindent 3. The variational approach to the geodesic equation\leaders \hbox to 1em{\hss .\hss }\hfill {\eightrm 7}\par 
\noindent 4. Some properties of the geodesics \leaders \hbox to 1em{\hss .\hss }\hfill {\eightrm 10}\par 

\head\totoc\nmb0{0}. Introduction \endhead 
If $M$ is a (not necessarily compact)
smooth finite dimensional manifold, the space $\M=\M(M)$
of all Riemannian metrics on $M$ can be endowed with
a structure of an infinite dimensional smooth manifold modeled
on the space $C^\infty_c(S^2T^*M)$ of symmetric $(0,2)$-tensor fields with
compact support. 
Analogously, the space $\Om^2_{\text{nd}}(M)$ of non degenerate 
2-forms on $M$, is an infinite dimensional smooth manifold, modeled 
on the space $\Om^2_c(M)$ of 2-forms with compact support.
See \cit!{6} and \cit!{7}. 

Here we consider the space of almost Hermitian structures on $M$, 
i.e\. the subset $\Cal H$ of $\M(M)\x \Om^2_{\text{nd}}(M)$ of those 
elements $(g,\om)$ such that the $(1,1)$-tensor field $J=g\i\om$ is an 
almost complex structure on $M$. 
The aim of this paper is to study the geometry of $\Cal H$. 
First we 
prove in section \nmb!{1} that $\Cal H$ is a splitting submanifold of 
the product $\M(M)\x \Om^2_{\text{nd}}(M)$. 
In section \nmb!{2}, 
after splitting the tangent space of the product in a form well 
adapted to our problem, we derive the equations that a curve in  
$\Cal H$ should satisfy in order to be a geodesic for the 
metric on $\Cal H$ induced by the product metric on 
$\M(M)\x \Om^2_{\text{nd}}(M)$. 
In section \nmb!{3} we give an 
independent variational derivation of the geodesic equations, by 
parameterizing elements of $\Cal H$ by automorphisms of $TM$, i.e\. by 
sections of the bundle $GL(TM,TM)$. 
Finally, section \nmb!{4} is devoted to the study of the geodesic equations
found in section \nmb!{2}. 
We were not able to find an explicit 
solution, nevertheless we can give some explicit properties of the 
geodesics, see \nmb!{4.1} and \nmb!{4.2}.
The subspaces $\Cal H_J$, $\Cal H_\om$, and $\Cal H_g$ of all almost 
hermitian structures with fixed almost complex structure $J$, 2-form 
$\om$, or metric $g$, respectively, are splitting submanifolds  of 
$\Cal H$. 
This follows from splitting the bundle. 
Some other 
interesting subsets (hermitian structures, K\"ahler structures, with 
symplectic $\om$) are much more difficult to treat: we do not know 
whether they are submanifolds, since the differential operators 
describing them are complicated in the charts we use. 
Moreover 
$\Cal H_J$ is totally geodesic in $\M(M)$, and the other two are 
totally geodesic in the manifold of all metrics with a fixed volume 
form. 
The geodesics are known explicitly in all three cases.

\head\totoc\nmb0{1}. Almost hermitian structures
\endhead
 
\subhead{\nmb.{1.1}. Almost hermitian structures}\endsubhead
Let $M$ be a smooth manifold of even dimension $n=2m$. 
Let $g$ be a
Riemannian metric on $M$ and let $\om$ be a non degenerate 2-form on
$M$; both of them will be regarded as fiber bilinear functionals on $TM$ or
as fiber linear isomorphisms $TM\to T^*M$ without any change of
notation. 
The symmetry of $g$ is then expressed $g^t=g$ where the 
transposed $g^t$ is given by $TM\to T^{**}M @>{g^*}>> T^*M$, and we 
also have $\om^t=-\om$
 
Then we consider the endomorphism $J:= g\i\om\colon TM\to TM$ which 
satisfies $J^*= \om^t(g^t)\i = -\om g\i$. 
Then the following conditions are equivalent:
\roster
\item $J$ is an isometry for $g$, i\.e\. $g(JX,JY)=g(X,Y)$ for all
       $X,Y\in T_xM$.
\item $g\i\om+\om\i g=0$.
\item $J$ is an almost complex structure, i\. e\.
       $J^2=- Id$ or equivalently $J\i = -J$.
\endroster
If these equivalent conditions are satisfied we say that $(g,\om)$ is
an \idx{\bf almost Hermitian structure}. 
 
{\bf Remark:} For a given almost complex structure $J$, condition 
\thetag1 is equivalent to the fact that $gJ$ is skew symmetric. 
So there is a bijective correspondence between the set of almost 
hermitian structures and the set of 
Riemannian almost complex structures.
 
\subhead{\nmb.{1.2}. The bundle of almost hermitian structures}\endsubhead
We consider the subspace
$$
\operatorname{Herm} :=\{(g,\om)\colon g\i\om+\om\i g=0\}\subset
     S^2_+T^*M\x \La^2_{\text{nd}}T^*M,
$$
where $S^2_+T^*M$ is the set of all positive definite symmetric 
2-tensors on $M$,
and we claim that it is a subbundle.
For that we consider the following commutative diagram:
$$\CD
\operatorname{Herm} @>pr_1>> S^2_+T^*M \\
@V{pr_2}VV      @VV{\pi}V\\
\La^2_{\text{nd}}T^*M @>>\pi> M.
\endCD\tag1$$
 
\proclaim{Lemma}
This is a double fiber bundle, where the standard
fiber of $pr_1$ is the homogeneous space $O(2m,\Bbb R)/U(m)$ and the
standard fiber of $pr_2$ is the homogeneous space
$Sp(m,\Bbb R)/U(m)$.
 
In fact all the fiber bundles of diagram \thetag1 are associated
bundles for the linear frame bundle $GL(\Bbb R^{2m},TM)$, where the
structure group $GL(2m)$ acts from the left on the typical fibers 
given in the diagram
$$\CD
GL(2m,\Bbb R)/U(m) @>>> GL(2m,\Bbb R)/O(2m,\Bbb R)\\
@VVV      @VVV\\
GL(2m,\Bbb R)/Sp(m,\Bbb R) @>>> \{Id\}.
\endCD\tag 2 $$
\endproclaim
 
\demo{Proof}
If we fix a metric $g$, then in an orthonormal frame of $TM|U$ for an 
open subset $U\subset M$, 
there is of course a 2-form $\om_0$ such that $(g,\om_0)$ is a local
section of $\operatorname{Herm}$. 
If $(g,\om_1)$ is another local section with the same
$g$ then there is a $g$-isometric local isomorphism $f$ of $TM|U$
with possibly smaller $U$
such that $f^*\om_1=\om_0$ which is fiberwise
unique up to multiplication from
the right by an element of $O(2m,\Bbb R)\cap Sp(m,\Bbb R) =U(m)$.
 
If we fix on the other hand a non degenerate 2-form $\om$, then there
is a frame of $TM|U$ such that $\om$ takes the usual
standard form of a symplectic structure (the frame can be chosen
holonomic if and only if $d\om=0$). 
Then obviously there is a metric
$g_0$ (constant in that frame) such that $(g_0,\om)$ is a local section
of $\operatorname{Herm}$. 
If $(g_1,\om)$ is another local section then there is a
fiberwise symplectic isomorphism $f$ of $TM|U$ (with possibly smaller
$U$) such that $f^*g_1 =g_0$, which is fiberwise unique up to right
multiplication by an element of $Sp(m,\Bbb R)\cap O(2m,\Bbb R)=U(m)$.

To check the last statement we consider the diagram
$$\CD
GL(\Bbb R^{2m},TM)\x GL(2m,\Bbb R)/O(2m,\Bbb R) @>A>> S^2_+T^*M\\
@VVV @| \\
\frac{GL(\Bbb R^{2m},TM)\x GL(2m,\Bbb R)/O(2m,\Bbb R)}{GL(2m,\Bbb R)}
	@>{\cong}>>  S^2_+T^*M
\endCD$$
where $A(u,g.O(2m,\Bbb R))(X,Y)=u\i(X)^tgg^tu\i(Y)$ for 
$u\in GL(\Bbb R^2m,T_xM)$ and $X, Y\in T_xM$.

Likewise we have 
$$\CD
GL(\Bbb R^{2m},TM)\x GL(2m,\Bbb R)/Sp(m,\Bbb R) @>B>> \La^2_{\text{nd}}T^*M\\
@VVV @| \\
\frac{GL(\Bbb R^{2m},TM)\x GL(2m,\Bbb R)/Sp(m,\Bbb R)}{GL(2m,\Bbb R)}
	@>{\cong}>>  \La^2_{\text{nd}}T^*M
\endCD$$
where $B(u,g.Sp(m,\Bbb R))(X,Y)=u\i(X)^tgJg^tu\i(Y)$.
\qed\enddemo

\subhead\nmb.{1.3}. The manifold of almost hermitian structures \endsubhead 
The space of almost hermitian structures on $M$ is just the
space $\Cal H:= C^\infty(\operatorname{Herm})$ of smooth sections of the fiber bundle $\operatorname{Herm}$.
Since $\operatorname{Herm}$ is a subbundle 
of $S^2_+T^*M\x \La^2_{\text{nd}}T^*M$ and since the latter is an 
open subbundle of the vector bundle $S^2T^*M\x \La^2T^*M$ we 
see that $\Cal H$ is a splitting submanifold of $\M(M)\x \Om^2_{\text{nd}}(M)$, 
by the following lemma. The splitting property will also follow 
directly in section \nmb!{2.1}.

\proclaim{Lemma} Let $(E,p,M,S)$ and $(E',p',M,S')$ be two fiber 
bundles over $M$ and let $i\colon E\to E'$ be a fiber respecting embedding. 
Then the following embeddings of spaces of sections are splitting 
smooth submanifolds:
$$
C^\infty(E) @>{i_*}>> C^\infty(E') \hookrightarrow C^\infty(M,E').
$$
\endproclaim

\demo{Proof}
This is a variant of the results 10.6 and 10.10 in \cit!{10} 
and the proof is similar to the ones given there, by direct 
finite dimensional construction.
\qed\enddemo
 
\subhead{\nmb.{1.4}. Metrics on $\Cal H= C^\infty(\operatorname{Herm})$}\endsubhead
On the space $\M(M)\x \Om^2_{\text{nd}}(M)$ there are many pseudo
Riemannian metrics which are invariant under the diffeomorphism group
and which are of first order.
We shall consider the product metric
$$\multline
G_{(g,\om)}((h_1,\ph_1),(h_2,\ph_2)) =\\
=\int_M \tr(g\i h_1 g\i h_2)\Vol(g) 
     +\int_M \tr(\om\i \ph_1\om\i \ph_2)\Vol(\om),
\endmultline$$
where 
$(h_1,\ph_1),(h_2,\ph_2)\in T_g\M(M)\x T_\om\Om^2_{\text{nd}}(M)= 
C^\infty_c(S^2T^*M)\x\Om^2_c(M)$.
Note that $\Vol(g)=\Vol(\om)$ if 
$(g,\om)\in \Cal H$ since $Sp(m,\Bbb R)$ and $O(2m,\Bbb R)$ are both 
subgroups of $SL(2m,\Bbb R)$. 
So the restriction of $G$ to $\Cal H$ 
is given by 
$$
G_{(g,\om)}((h_1,\ph_1),(h_2,\ph_2)) :=
\int_M\Bigl( \tr(g\i h_1 g\i h_2) + \tr(\om\i \ph_1
\om\i \ph_2)\Bigr)\Vol(g).
$$

\head\totoc\nmb0{2}. The geodesic equation in $\Cal H$ \endhead 

\subhead{\nmb.{2.1}. Splitting the tangent space of the space of almost
hermitian structures}\endsubhead
Let
$(g,\om)\in \Cal H=C^\infty(\operatorname{Herm})\subset 
     \M(M)\x \Om^2_{\text{nd}}(M)$
so that $F(g,\om):= g\i\om + \om\i g=0$.
Then for $(h,\al)\in C^\infty_c(S^2T^*M\x \La^2T^*M)$ we put
$H:=g\i h$, $A := \om\i \al$, and of course $J=g\i\om$.
We have $(h,\al)\in T_{(g,\om)}\Cal H$ if and only if
$$
dF(g,\om)(h,\al)=-g\i h g\i \om + g\i\al + \om\i h
     - \om\i\al\om\i g = 0, 
$$
which is easily seen to be equivalent to $JH+HJ=JA+AJ$. 
Note that this implies $\tr(H)=\tr(A)$.
 
For $(g,\om)\in\Cal H$ we have the following $G$-orthogonal decomposition
of the tangent space to $\M(M)\x \Om^2_{\text{nd}}(M)$ at $(g,\om)$:
$$\align
T_{(g,\om)}(\M(M)&\x \Om^2_{\text{nd}}(M))
     = \Cal N^1_{(g,\om)} \oplus \Cal N^2_{(g,\om)}
     \oplus \Cal N^3_{(g,\om)} \oplus \Cal N^4_{(g,\om)} \\
\Cal N^1_{(g,\om)} &= \{(H,0)\colon JHJ=H\}\\
\Cal N^2_{(g,\om)} &= \{(0,A)\colon JAJ=A\}\\
\Cal N^3_{(g,\om)} &= \{(H,A)\colon JHJ=-H\text{ and }H=A\}\\
\Cal N^4_{(g,\om)} &= \{(H,A)\colon JHJ=-H\text{ and }H=-A\}
\endalign$$
The tangent space to $\Cal H$ is given by
$T_{(g,\om)}\Cal H = \Cal N^1_{(g,\om)} \oplus \Cal N^2_{(g,\om)}
     \oplus \Cal N^3_{(g,\om)}$
and its $G$-orthogonal complement is given by
$(T_{(g,\om)}\Cal H)^\bot = \Cal N^4_{(g,\om)}$.
The restriction of the pseudometric to $\Cal H$ is then non 
degenerate. 

The projectors on these subspaces can easily be constructed and in
particular we have
the orthogonal projectors from
$T_{(g,\om)}\M(M)\x \Om^2_{\text{nd}}(M)$ to the tangent space
$T_{(g,\om)}\Cal H$ and to the $G$-orthogonal complement
$$\gather
Pr^{\text{T}}_{(g,\om)}(H,A):=
     \left(\frac{3H+JHJ+A-JAJ}4,\frac{3A+JAJ+H-JHJ}4\right)\\
Pr^{\bot}_{(g,\om)}(H,A):=
     \left(\frac{H-JHJ-A+JAJ}4,\frac{A-JAJ-H+JHJ}4\right).
\endgather$$

\subhead\nmb.{2.2}. The geodesic equation \endsubhead
The space $\Cal H$ is a splitting submanifold of 
$\M(M)\x \Om^2_{\text{nd}}(M)$ and the tangent space splits nicely in 
the direct sum of the tangential and the orthogonal part, by 
\nmb!{2.1}; note that the projection operators are algebraic. 
The 
tangential projection of the covariant derivative $\nabla_\xi\et$ of 
smooth vector fields on $\M(M)\x \Om^2_{\text{nd}}(M)$ which along 
$\Cal H$ are tangential to $\Cal H$, is thus again a smooth vector 
field, and the usual proof of Gau\ss' formula involving the six 
term expression of the Levi-Civita covariant derivative shows that 
$Pr^T\nabla_\xi\et$ is the smooth Levi-Civita covariant derivative of 
$G$ on $\Cal H$. 

So a curve $\si(t)=(g(t),\om(t))$ is a geodesic for the 
induced metric if and only if the covariant derivative of its tangent 
vector $\si_t$ in $\M(M)\x \Om^2_{\text{nd}}(M)$ is everywhere 
orthogonal to $\Cal H$; i\. e\. 
$\nabla_{\partial_t}\si_t\in\Cal N_\si^4$, or equivalently 
$Pr^T_\si(\nabla_{\partial_t}\si_t)= 0$.

If we put $X:=g\i g_t$ and $W:=\om\i\om_t$ and use the Christoffel 
form from
\cit!{7}, (2.3 for $\al=1/n$), these 
conditions become, respectively:
$$\cases 
& JX_t + \tfrac12 \tr(X)JX = X_tJ + \tfrac12 \tr(X)XJ\\
& X_t + W_t = - \tfrac12 \tr(X)X - \tfrac12 \tr(W)W + 
	\tfrac14(\tr(X^2)+\tr(W^2))Id
\endcases\tag{1}$$

$$\cases
&3X_t + \tfrac32 \tr(X)X - \tfrac12 \tr(X^2)Id + JX_tJ 
	+ \tfrac12 \tr(X)JXJ +\\
&+ W_t + \tfrac12 \tr(W)W - \tfrac12 \tr(W^2)Id - JW_tJ 
	- \tfrac12 \tr(W)JWJ = 0\\
&3W_t + \tfrac32 \tr(W)W - \tfrac12 \tr(W^2)Id + JW_tJ 
	+ \tfrac12 \tr(W)JWJ +\\
&+ X_t + \tfrac12 \tr(X)X - \tfrac12 \tr(X^2)Id - JX_tJ 
	- \tfrac12 \tr(X)JXJ = 0
\endcases\tag{2}$$

\subhead\nmb.{2.3}. The submanifold $\Cal H_J$ of $\Cal H$ \endsubhead
For a fixed almost complex structure $J$ on $M$ let us consider
$\Cal H_J := \{(g,\om)\in \Cal H\colon  g\i\om=J\}$, the space of
almost hermitian structures with almost complex structure
$J$. 
The tangent space is $T_{(g,\om)}\Cal H_J = 
\Cal N^3_{(g,\om)}$, see \nmb!{2.1}

Via the first projection the space $\Cal H_J$ is diffeomorphic to the 
submanifold of $\M(M)$ consisting of all $g$ making $J$ an isometry. 
It is a totally geodesic submanifold in the 
Riemannian manifold $(\M(M),G)$. 
This follows from the general 
result: 
 
\proclaim{Lemma} Let $A\colon TM\to TM$ be a vector bundle
isomorphism covering the identity. 
Then the space of all Riemannian
metrics $g$ on $M$ such that $g(AX,AY)=g(X,Y)$ for all $X,Y\in T_xM$
is a geodesically closed submanifold of $(\M(M),G)$. 
This is also
true for bilinear structures, or metrics with fixed signature.
\endproclaim
 
\demo{Proof}
The space of these $g$ is the space of sections of the open subbundle 
$S^2_+T^*M\cap\{g\colon  A^*\o g\o A= g\}$ of the obvious vector
subbundle, which is a smooth manifold. 
A tangent vector $h$ at such $g$
is a tensor field with compact support with $A^*\o h\o A= h$; for
$H=g\i h$ this is equivalent to $A\i\o H\o A= H$.
By \cit!{6},~3.2, the geodesic in $\M(M)$ starting at $g$ in the 
direction $h$ is the curve 
$$
g(t)=g\;e^{a(t)Id+b(t)H_0},
$$
where $H_0$ is the traceless part of $H$ and $a(t)$ and $b(t)$ are 
real valued functions. 
Then if $g$ and $h$ are $A$-invariant, so is 
the whole geodesic.  
\qed\enddemo

\subhead\nmb.{2.4}. The submanifold $\Cal H_\om$ \endsubhead
For a fixed non degenerate 2-form $\om$ we consider $\Cal H_\om := 
\M(M)\x \{\om\}\cap\Cal H$, the space of almost
hermitian structures with fixed $\om$. 
It
is the space of sections of the pullback bundle
$\om^*(\operatorname{Herm},pr_2,\La^2_{\text{nd}}T^*M)$ in terms of 
\nmb!{1.2},
so a smooth manifold.
Its tangent space is $T_g\Cal H_\om = \Cal N^1_{(g,\om)}$.
The space $\Cal H_\om$ is a submanifold $\M_{\Vol(\om)}$  of the space 
of all metrics with fixed volume equal to $\Vol(\om)$, see 
\nmb!{1.4}. 

$\Cal H_\om$ is a geodesically closed submanifold of $\M_{\Vol(\om)}$, 
see Blair \cit!{1}, \cit!{2}. 
The geodesics of $\M_{\Vol(\om)}$ have been 
determined by Ebin \cit!{3}, see also \cit!{4}. 
The geodesic starting at $g$ in the direction $h$ is given by 
$$
g(t)= g\;e^{tH_0}.
$$
Note that $\M_{\Vol(\om)}$ is not geodesically closed in $\M$. 
These results are given for 
compact $M$, but clearly they continue to hold for noncompact $M$ in 
our setting.

\subhead\nmb.{2.5}. The submanifold $\Cal H_g$ \endsubhead
For a fixed metric $g$ we may consider 
$\Cal H_g := \{g\}\x \Om^2_{\text{nd}}(M)\cap \Cal H$,
the space of almost hermitian structures with fixed $g$. 
It is
the space of all sections of the pullback bundle
$g^* (H,pr_1,S^2_+T^*M)$ in the notation of \nmb!{1.2}, so it is a
smooth manifold. 
For the tangent space it easily follows that
$T_\om\Cal H_g = \Cal N^2_{(g,\om)}$.
Since the situation is symmetric with respect to $g$ or $\om$, it 
follows from \nmb!{2.4} that $\Cal H_g$ is a geodesically closed submanifold of the 
space $\Om^2_{\text{nd}}(M)_{\Vol(g)}$ of all almost symplectic 
structures with fixed volume equal to $\Vol(g)$, and the geodesics in 
$\Om^2_{\text{nd}}(M)_{\Vol(g)}$ are given by 
$$
\om(t)=\om\; e^{tA_0}.
$$

\head\totoc\nmb0{3}. The variational approach to the geodesic 
equation\endhead

\subhead{\nmb.{3.1}}\endsubhead 
It is not so easy to find an adapted chart 
for the subbundle $\operatorname{Herm}\subset S^2_+T^*M\x \La^2_{\text{nd}}T^*M$ which 
would allow us to parameterize curves and their variations in 
$\Cal H$. 
In order to achieve this parameterization we will use the 
following scheme.

We consider the bundle $GL(TM,TM)$	of all isomorphisms of the tangent 
bundle. 
It is an open submanifold of the vector bundle $L(TM,TM)$. 
Any (fixed) almost hermitian structure $(g_0,\om_0)\in\Cal H$ induces 
a smooth mapping 
$$\gather
\ph=\ph_{(g_0,\om_0)}\colon  GL(TM,TM) \to \operatorname{Herm} \subset 
	S^2_+T^*M\x \La^2_{\text{nd}}T^*M \\
\ph(f)=\ph_{(g_0,\om_0)}(f)= (f^*g_0f,f^*\om_0f).
\endgather$$
The corresponding push forward mapping between the spaces of sections 
will be denoted by 
$$\gather
\Ph\colon \Cal G:= C^\infty(GL(TM,TM)) 
	\to \Cal H=C^\infty(\operatorname{Herm})\subset \Cal M\x \Om^2_{\text{nd}}(M),\\
\Ph(f)=\Ph_{(g_0,\om_0)}(f)= \ph\o f = (f^*g_0f,f^*\om_0f).
\endgather$$
 
\remark{Remark} The mapping 
$\Ph\colon \Cal G=C^\infty(GL(TM,TM))\to \Cal H$ is {\bf not} surjective in 
general. 
The mapping $\ph\colon GL(TM,TM)\to \operatorname{Herm}$ is the projection of a 
fiber bundle with typical fiber $U(m)$. 
Since $U(m)$ has nontrivial 
homotopy trying to lift a section $s\colon M\to \operatorname{Herm}$ over $\ph$ will meet 
obstructions in general.
\endremark

\proclaim{\nmb.{3.2}. Lemma} For every curve $(g(t),\om(t))$ in 
$\Cal H$ and also for every variation $(g(t,s),\om(t,s))$ of such a curve 
in $\Cal H$ with $(g(0),\om(0))=(g_0,\om_0)$ there is a curve $f(t)$ or 
variation $f(t,s)$ in $\Cal G= C^\infty(GL(TM,TM))$ with 
$(g,\om)=\Ph_{(g_0,\om_0)}(f)$.
\endproclaim

\demo{Proof}
As noticed above $\ph\colon GL(TM,TM)\to \operatorname{Herm}$ is the projection of a smooth 
fiber bundle with compact fiber type $U(m)$. 
We choose a generalized 
connection for this bundle, see \cit!{11} or 
\cit!{9}, section~9. 
Its parallel transport 
$$
\operatorname{Pt}(c,t)\colon GL(TM,TM)_{c(0)}\to GL(TM,TM)_{c(t)}
$$ 
is globally defined for each curve 
$c\colon \Bbb R\to \operatorname{Herm}$, and it is smooth in 
the choice of the curve, see loc\. cit\.

Then we just define 
$f(t):= \operatorname{Pt}((g(x,\quad),\om(x,\quad)),t)
	\operatorname{Id}_{T_xM}$
and $f(\quad,t)$ will be a curve in $\Cal G= C^\infty(GL(TM,TM))$ 
with $(g(t),\om(t))=\Ph_{(g_0,\om_0)}(f(t))$, and $f(x,t)$ varies in 
$t$ only for those $x$ where also $g(x,t)$ or $\om(x,t)$ varies. 
So 
the compact support condition required of smooth curves of sections 
is automatically satisfied. 

For a variation of a curve we first define in turn
$$\align
f(x,t,0):&= \operatorname{Pt}((g(x,\quad,0),\om(x,\quad,0)),t)
	\operatorname{Id}_{T_xM},\\
f(x,t,s) :&= 
\operatorname{Pt}((g(x,t,\quad),\om(x,t,\quad)),s)f(x,t,0).\qed
\endalign$$
\enddemo
 
\subhead{\nmb.{3.3}}\endsubhead 
Let $(g(t),\om(t))$ be a smooth curve in 
$\Cal H=C^\infty(\operatorname{Herm})$, 
so it is smooth $M\x \Bbb R\to \operatorname{Herm}$ and for each 
compact $[a,b]\subset \Bbb R$ there is a compact set $K\subset M$ 
such that $(g(x,t),\om(x,t))$ is constant in $t\in[a,b]$ for each 
$x\in M\setminus K$, see (\cit!{11},~6.2, a slight mistake there). 
Then its energy with respect to the 
metric $G$ of \nmb!{1.4} is given by 
$$
E_a^b(g,\om) = \int_a^b\int_M \left(\tr(g\i g_t g\i g_t) + 
	\tr(\om\i \om_t \om\i \om_t)\right) \Vol(g) dt. \tag1
$$
Since we cannot parameterize curves and their variations in $\Cal H$ 
explicitly we will parameterize them with the help of 
$\Ph=\Ph_{(g_0,\om_0)}\colon \Cal G\to \Cal H$, where 
$(g_0,\om_0)=(g(0),\om(0))$. 
So let $f(t)$ be a smooth curve in 
$\Cal G= C^\infty(GL(TM,TM))$. 
Then we have 
$$\gather
\Ph(f)=\Ph_{(g_0,\om_0)}(f)=(f^*g_0f,f^*\om_0f),\\
T_f\Ph(f_t)=(f^*_tg_0f+f^*g_0f_t,f^*_t\om_0f+f^*\om_0f_t),
\endgather$$
so the energy of the curve $\Ph(f(t))$ in $\Cal H$ is given by
$$\align
\text{\thetag2}&\qquad E_a^b(\Ph(f)) =\\ 
&= \int_a^b\int_M \Bigl(
	\tr(f\i g_0\i (f^*)\i(f^*_tg_0f+f^*g_0f_t)
		f\i g_0\i (f^*)\i(f^*_tg_0f+f^*g_0f_t)) \\
&\quad+\tr(f\i \om_0\i (f^*)\i(f^*_t\om_0f+f^*\om_0f_t)
	f\i \om_0\i (f^*)\i(f^*_t\om_0f+f^*g_0f_t)) 
     \Bigr)\\
&\qquad\det(f)\Vol(g_0)dt.
\endalign$$

\proclaim{Lemma} A curve $f(t)$ in $\Cal G$ is a critical point of 
the functional \thetag2 if and only if $\Ph(f(t))$ is a critical 
point of the functional \thetag1, i\. e\. $t\mapsto \Ph(f(t))$ is a 
geodesic.
\endproclaim

\demo{Proof}
An (infinitesimal) variation in $\Cal G$ can (with the help of a 
connection for $\ph\colon GL(TM,TM)\to \operatorname{Herm}$) be written as a sum of two 
variations: the horizontal one corresponds exactly to a variation of 
$\Ph(f)$ in $\Cal H$, and along the vertical one the functional is
stationery anyhow.
\qed\enddemo

\proclaim{\nmb.{3.4}. Lemma} In the setting of \nmb!{3.3}, for a 
variation $f(t,s)\in\Cal G$ with fixed endpoints we have the first 
`variation formula' 
$$\align
\left.\frac{\partial}{\partial s}\right|_0 
	&E_a^b(\Ph_{(g_0,\om_0)}(f(\quad,s))) = \\
&= \int_a^b\int_M \tr\Bigl( E(g_0,\om_0,f;t)
	f_sf\i\Bigr)\det(f)\Vol(g)\,dt.
\endalign$$
where 
$$\align
E(g_0,\om_0,f;t) &=	- 2f_{tt}f\i - 2g_0\i(f^*)\i f^*_{tt}g_0 
	+ 2f_tf\i f_tf\i \\
&\quad - 2\tr(f_tf\i)f_tf\i + \tr(f_tf\i f_tf\i)Id  
	- 2g_0\i (f^*)\i f^*_tg_0f_tf\i \\
&\quad + 2g_0\i (f^*)\i f^*_t(f^*)\i f^*_tg_0 
	+ 2f_tf\i g_0\i (f^*)\i f^*_tg_0 \\
&\quad - 2\tr(f_tf\i)g_0\i (f^*)\i f^*_tg_0 
	+ \tr(g_0\i (f^*)\i f^*_tg_0f_tf\i)Id \\
&\quad - 2f_{tt}f\i - 2\om_0\i(f^*)\i f^*_{tt}\om_0 
	+ 2f_tf\i f_tf\i \\
&\quad - 2\tr(f_tf\i)f_tf\i + \tr(f_tf\i f_tf\i)Id  \\
&\quad - 2\om_0\i (f^*)\i f^*_t\om_0f_tf\i 
	+ 2\om_0\i (f^*)\i f^*_t(f^*)\i f^*_t\om_0 \\
&\quad + 2f_tf\i \om_0\i (f^*)\i f^*_t\om_0 
	- 2\tr(f_tf\i)\om_0\i (f^*)\i f^*_t\om_0 \\
&\quad + \tr(\om_0\i (f^*)\i f^*_t\om_0f_tf\i)Id
\endalign$$
\endproclaim

\demo{Proof}
This is a long but straightforward computation.
We may interchange $\frac{\partial}{\partial s}|_0$ with the
first integral since this is finite dimensional analysis, and we
may interchange it with the second one, since $\int_M$ is a
continuous linear functional on the space of all smooth
densities with compact support on $M$, by the chain rule.
Then we use that $\tr_*$ is linear and continuous,
$d(\Vol)(g)h= \frac12\tr(g\i h)\Vol(g)$, and that
$d((\quad)\i)_*(g)h= -g\i h g\i$ and partial integration; there are 
no boundary terms since we assumed the variation to have fixed 
endpoints.
\qed\enddemo

\proclaim{\nmb.{3.5}. Lemma} For curve $f(t)$ in $\Cal G$ the curve 
$\Ph_{(g_0,\om_0)}(f(t))$ is a geodesic in $(\Cal H, G)$ if and only 
if $f(t)$ satisfies the following equation:
$$
E(g_0,\om_0,f,t) = 0.
$$
\endproclaim

\demo{Proof}
This follows from \nmb!{3.4} since the integral in \thetag3 
describes a nondegenerate inner product on $\Cal G$, given by
$$
G_f(h,k)=\int_M\tr(hf\i kf\i)\det(f)\Vol(g_0).\qed
$$
\enddemo

\subhead\nmb.{3.6}. Comparison with section \nmb!{2} \endsubhead
Let $\Ph_{(g_0,\om_0)}(f(t))=(f^*g_0f,f^*\om_0f)=:(g(t),\om(t))$. 
Then the expressions used in section \nmb!{2} become 
$$\align
X&= g\i g_t =f\i g_0\i (f^*)\i f^*_t g_0 f + f\i f_t,\\
W&= \om\i \om_t =f\i \om_0\i (f^*)\i f^*_t \om_0 f + f\i f_t,\\
J&= g\i \om = f\i g_0\i \om_0 f.
\endalign$$
If we compute $X_t$, $W_t$ and insert this into the second equation 
of \nmb!{2.2}.(1), we get exactly $f\i E(g_0,\om_0,f,t) f =0$.
So we get the same geodesic equation as in \nmb!{2.2}

\head\totoc\nmb0{4}. Some properties of the geodesics \endhead

We are not able to give the explicit solution of the geodesic 
equation on $\Cal H$. 
But we can give some explicit formulas of the 
time evolution of some functions of the structures.

\proclaim{\nmb.{4.1}. Proposition}
Let $(g(t),\om(t))$ be the geodesic of $\Cal H$ starting at 
$(g_0,\om_0)$ in the direction $(h,\al)$, let $H=g_0\i h$ and 
$A=\om_0\i\al$, and let $(X,W)=(g\i g_t, \om\i\om_t)$ as in 
\nmb!{2.2}. 
Then we have 
$$
\tr(X^2)+\tr(W^2) = (\det(g_0\i g))^{-1/2}(\tr(H^2)+\tr(A^2)).
$$
\endproclaim

\demo{Proof}
Since $(g,\om)$ is in $\Cal H$ its tangent vector $(g_t,\om_t)$ 
satisfies 
$$\gather
JXJ - X = JWJ - W \tag1 \\
\tr(X) = \tr(W) \tag 2
\endgather$$
so that \nmb!{2.2}.\thetag 2 becomes
$$\cases
&3X_t + JX_tJ + W_t - JW_tJ +2\tr(X)X - \tfrac12(\tr(X^2)+\tr(W^2))Id = 0 \\
&3W_t + JW_tJ + X_t - JX_tJ +2\tr(W)W - \tfrac12(\tr(W^2)+\tr(X^2))Id = 0 
\endcases$$
We multiply now the first equation by $X$, the second 
equation by $W$ and add them to obtain 
$$
2(\tr(X_tX)+\tr(W_tW))+\tfrac12(\tr(X^2)+\tr(W^2))\tr(X)=0.
$$
 From that it is easy to see that the derivative of 
$(\tr(X^2)+\tr(W^2))(\det(g_0\i g))^{1/2}$ is zero, where we also use
$((\det(g_0\i g))^{1/2})_t = \tfrac12 (\det(g_0\i g))^{1/2}\tr(X)$.
\qed\enddemo

\proclaim{\nmb.{4.2}. Proposition}
Let $(g(t),\om(t))$ be the geodesic of $\Cal H$ starting at 
$(g_0,\om_0)$ in the direction $(h,\al)$, let $H=g_0\i h$ and 
$A=\om_0\i\al$, and let $(X,W)=(g\i g_t, \om\i\om_t)$. 
We put $p(t):= \tfrac12 (\det(g_0\i g))^{1/2}$. 
Then we have
$$\align
p(t) &= \tfrac n{32}(\tr(H^2)+\tr(A^2))t^2 + \tfrac12\tr(H)t +1, \tag 1\\
\tr(X) &= \tr(W) = 2\frac {p'(t)}{p(t)}\tag2 \\
&= 4 \frac{n(\tr(H^2)+\tr(A^2))t + 8\tr(H)}
	{n(\tr(H^2)+\tr(A^2))t^2 + 8\tr(H)t +32},\\
X+W &= \frac1{p(t)}\left(\frac14(\tr(H^2)+\tr(A^2))t Id + H + 
	A\right) \tag3
\endalign$$
\endproclaim

\demo{Proof}
We take the trace in the second expression of \nmb!{2.2}.\thetag 1 and use 
\thetag 2 from the proof of \nmb!{4.1} to obtain
$$
2\tr(X)' = - \tr(X)^2 + \tfrac n4(\tr(X^2)+\tr(W^2)).
$$
Inserting \nmb!{4.1} we get
$$
2\tr(X)' = - \tr(X)^2 + \tfrac n4\frac C{p(t)},
$$
where $C=\tr(H^2)+\tr(A^2).$
 From the proof of \nmb!{4.1} we have in turn
$$\align
p'(t) &= \tfrac12 p(t)\tr(X)\\
p''(t) &=	\tfrac 14 p(t)\tr(X)^2 + \tfrac12 p(t)\tr(X)' = 
\frac {nC}{16}.
\endalign$$
For the initial conditions $p(0)=1$ and $p'(0)=\tfrac12 \tr(H)$ this 
gives
$$
p(t) = \tfrac {nC}{32}t^2 + \tfrac12\tr(H)t +1
$$
and consequently assertions \thetag1 and \thetag2.

Now we take the tracefree part of the second expression 
in \nmb!{2.2}.\thetag 1
$$
(X+W)_0'= -\tfrac14 \tr(X+W)(X+W)_0,
$$
and by an argument similar to that used in the proof of
\cit!{7},~2.5 we get  
$$
X+W = a(t)Id + b(t)(H_0+A_0),\text{ where }
     a(t)= \frac {2\tr(X)}n = \tfrac 4n\frac{p'(t)}{p(t)},
$$ 
and it just remains to find $b(t)$ when $(H_0+A_0)\ne0$.

We have $(X+W)_0=b(t)(H_0+A_0)$, thus $(X+W)_0'=b'(t)(H_0+A_0)$. 
But we also know that $(X+W)_0'=-\tfrac n4 a(t)b(t)(H_0+A_0)$ 
and so we get
$$
\frac{b'}b = -\frac {na}4 = -\frac{p'}p
$$
with $p(0)=1$ and $b(0)=1$,
so $b=1/p$ and we get assertion \thetag3.
\qed\enddemo

\enddocument